\documentclass[a4paper,12pt,portrait]{article}
\usepackage{amsmath}
\usepackage{amssymb}
\usepackage{graphicx}
\usepackage{caption}
\usepackage{subcaption}
\usepackage[dvipsnames]{xcolor}
\makeatletter \oddsidemargin  -.1in \evensidemargin -.1in
\newcommand{\doublespacing}{\let\CS=\@currsize\renewcommand
	{\baselinestretch}{1.5}\tiny\CS}
\begin{document}
	\newcommand{\bea}{\begin{eqnarray}}
		\newcommand{\eea}{\end{eqnarray}}
	\newcommand{\nn}{\nonumber}	
	\newcommand{\bee}{\begin{eqnarray*}}
		\newcommand{\eee}{\end{eqnarray*}}
	\newcommand{\lb}{\label}
	\newcommand{\nii}{\noindent}
	\newcommand{\ii}{\indent}
	\newtheorem{thm}{Theorem}[section]
	\newtheorem{cor}{Corollary}[section]
	\newtheorem{lem}{Lemma}[section]
	\newtheorem{rem}{Remark}[section]
	\renewcommand{\theequation}{\thesection.\arabic{equation}}
	\title{\bf Inadmissibility Results for the Selected Hazard Rates  }
	\author{Brijesh Kumar Jha$^a$\thanks {Email address: jha.brijesh5@gmail.com}, Ajaya Kumar Mahapatra$^b$\thanks {Email address (corresponding author):
			ajayaiter@gmail.com,~ajayamohapatra@soa.ac.in}~ and Suchandan Kayal$^c$\thanks {Email address: kayals@nitrkl.ac.in,~suchandan.kayal@gmail.com}}
	\date{}
	\maketitle
	\begin{quote}
		$a.$ {\it \small  Department of Mathematics, Siksha `O' Anusandhan University, Bhubaneswar-751030, India}\\
		$b.$ {\it \small  Centre for Applied Mathematics and Computing, Siksha `O' Anusandhan University, Bhubaneswar-751030, India}\\
		$c.$ {\it \small  Department of Mathematics, National Institute of
			Technology Rourkela, Rourkela-769008, India}
	\end{quote}
	
	\begin{center}
		\textbf{\textit{Final version will be published in ``Statistics''}}
	\end{center}
	 
	\begin{abstract}
		Let us consider $k ~(\ge 2)$ independent populations $\Pi_1, \ldots,\Pi_k$, where $\Pi_i$ follows exponential distribution with hazard rate ${\sigma_i},$ ($i = 1,\ldots,k$). Suppose
		$Y_{i1},\ldots, Y_{in}$ be a random sample
		of size $n$ drawn from the $i$th population $\Pi_i$, where
		$i = 1,\ldots,k.$ For $i = 1,\ldots,k$, consider
			$Y_i=\sum_{j=1}^nY_{ij}$. The natural selection
		rule is to select a population associated with
		the largest sample mean. That is, $\Pi_i$, ($i = 1,\ldots,k$)  is
		selected if $Y_i=\max(Y_1,\ldots,Y_k)$. Based on this selection rule,
		a population is chosen. Then, we consider the estimation of the hazard 
		rate of the selected population with respect to the entropy loss function. Some natural
		estimators are proposed. The minimaxity of a
		natural estimator is established. Improved estimators improving upon the natural estimators are derived. Finally, numerical study is carried out in order to compare the proposed estimators in terms of the risk values.\vspace{2mm}
		
	\end{abstract}
	
	\section{Introduction}
	In estimation after selection problem, a population
	is selected according to some selection procedures. Then, it is of interest to investigate the unknown parameter(s) or characteristic(s) of the selected population. For optimality of the selection procedures, one can refer to Bahadur and Goodman \cite{B3},
	Lehmann \cite{B9} and Eaton \cite{B6}. 
	Note that the parameters of interest of the selected population are usually the mean, reliability, hazard rate etc.
	In the literature, these problems are usually dubbed as the
	estimation after selection problems.  There are a lot of applications of the selected estimation problems. A leading smartphone brand may choose processor chipsets within the same price range from processor manufacturing companies like Qualcomm, Mediatek etc. on the basis of clock speed and then may estimate the hazard rate of the chipset from the selected population. Please note that clock speed decides the number of instructions that can be processed by a processor in one second. A company may be interested to install LED bulbs in his office from a selected brand on the basis of power consumption and then may go for estimation of failure rate of the bulbs.  The concept of hazard rate plays an important role in designing the safe systems. Hazard rate describes the instantaneous risk of observing the event of interest over time. Here, instantaneous risk stands for instantaneous failure.
	In this paper, we consider the estimation of the hazard rate from a selected exponential population with respect to the entropy loss function. It may be noted that the hazard rate of the exponential distribution is the reciprocal of the scale parameter.\\
	
	Let $\Pi_1,\ldots,\Pi_k$ be $k~(\ge 2)$ independent populations, where $\Pi_i$ has exponential distribution with probability density function given by 
	\begin{eqnarray*}
		f(y:\sigma_i)=\sigma_i~e^{-\sigma_i y},~y>0,~\sigma_i>0,~ i = 1,\ldots,k.
	\end{eqnarray*}
	Here, $\sigma_i$ is the hazard rate of the $i$th exponential population $\Pi_i$.
	Consider a random sample $Y_{i1},\ldots, Y_{in}$ taken from the $i$th population $\Pi_i, i = 1,\ldots,k .$
	Define $Y_i=\sum_{j=1}^n Y_{ij}$.
	It can be seen that $(Y_1,\ldots,Y_k)$ are independent gamma distributed random variables with a common shape parameter $n$ and different scale parameters $(\sigma_1^{-1},\ldots,\sigma_k^{-1})$ respectively with density, given by 
	\begin{eqnarray}
		g(y:\sigma_i,n)=\frac{\sigma_i^n}{\Gamma(n)}~y^{n-1}e^{-\sigma_i
			y},~y>0,~\sigma_i>0,~i = 1,\ldots,k.
		\label{2}
	\end{eqnarray}
	Note that if $i$ is omitted from Equation (\ref{2}), we obtain the probability density function of Gamma($\sigma,n$) distribution with scale parameter $\sigma^{-1}$ and shape parameter $n$. Henceforth, we use this notation throughout the paper. Denote $\underline{Y} = (Y_1,\ldots,Y_k).$ Then, $\underline{Y}$ is a joint complete and sufficient statistic for estimating $\underline{\sigma} = (\sigma_1, \ldots, \sigma_k).$
	Define $Y_J = \max(Y_1,\ldots,Y_k)$, and let $\sigma_{J}$ denote the hazard rate associated with $Y_J.$ The natural selection rule considered here is to select a population with the  largest sample mean. That is, 
	if $Y_i=\max(Y_1,\ldots,Y_k)$, then we select the population  $\Pi_i.$
	Thus, the hazard rate of the selected population  in the given problem is
	\begin{eqnarray}
		\sigma_{J}=\sum_{i=1}^k \sigma_i I_i,
		\label{ch5neh1}\end{eqnarray}
	where $I_i$ is the usual indicator function, given by 	$I_i =1,~\mbox{if}~ Y_i=Y_{J};$ $I_i =0$, otherwise. 
	Without loss of generality, we ignore the case of ties, as this has
	zero probability. Let us consider entropy loss function as
	\begin{eqnarray}
		L(\underline{\sigma}, 
		d)=\frac{d}{\sigma_{J}} - \ln\left(\frac{d}{\sigma_{J}}\right) - 1,
		\label{ch5.1.1}\end{eqnarray} where $d$ is an
	estimator of $\sigma_{J}$ and $\ln$ denotes the natural logarithm with base $e$. Suppose a basic part is selected from $k$ brands in the market and is used as an integral unit of a larger system. Then, it is of priority to have a better estimate of the selected part in order to acquire a precise estimate of the hazard rate of the entire system. It is worth to mention that the entropy loss function is simply an asymmetric loss that penalizes over estimates more than under estimates.\\
	\par
	The problem of estimation of scale parameter of a selected gamma population has been widely discussed in the literature. Some of the major contributions in this area were provided by  Sharma \cite{B16} and Vellaisamy \cite{B17, B18}. Vellaisamy \cite{B17, B18} considered the estimation of the scale parameter of the selected population. Using the method of differential inequalities, some improved estimators were derived. Misra et al. \cite{B11} studied the estimation of scale parameter of a selected gamma population by using a scale invariant squared error loss function. Sufficient conditions for inadmissibility of the scale equivariant estimators were established by the authors. Later, Misra et al. \cite{B12} developed some inadmissibility results for the scale
	parameters of selected gamma populations. Parsian and Nematollahi \cite{B14} studied estimation of scale parameter under the entropy loss function with restrictions
	to the principles of invariance and risk unbiasedness. An explicit form of the minimum risk scale-equivariant estimator under entropy loss function was obtained. Shariati and Nematollahi \cite{B13} considered the same problem under an asymmetric loss function  given by (\ref{ch5.1.1}). They derived the class of linear admissible estimators and showed that the UMRU (uniformly minimum risk unbiased) estimator is admissible.
	Arshad et al. \cite{B2} considered estimation of the mean of a selected gamma population under a scale invariant loss function. The problem of estimation of the scale parameter
	of a selected uniform population was addressed by Arshad et al. \cite{B1} with respect to the entropy loss function. They derived the UMRU estimator of scale parameter of the selected population. Two natural estimators of the scale parameters were proposed. In addition, they derived a general result for improving a scale invariant
	estimator of the scale parameter under the entropy loss function. It is seen that not much has been reported on the estimation of reliability or hazard rate of an exponential population. Kumar et al. \cite{B8} considered the problem of estimation of the reliability function from a selected subset of exponential populations. A class of improved estimators was obtained by solving a differential inequality which improves upon the natural estimators. Mahapatra et al. \cite{B10} studied the simultaneous estimation of the hazard rates of several exponential populations. Some new estimators were derived and their inadmissibility was proved. Some dominating estimators for the modified best scale equivariant estimator (BSEE) were also derived. Very recently, Jha et al. \cite{B7} considered the problem of estimation of the hazard rate of a selected exponential population and derived some improved estimators for the proposed natural estimators with respect to the quadratic loss function. However, the results obtained in the present paper are quite different from Jha et al. \cite{B7}.\\ 
	
	The paper is organized as follows. In Section $2$, Brewster-Zidek technique is used to derive a subclass of admissible estimators within a class of estimators. Minimaxity of a natural estimator is established in Section 3. An estimator improving upon the natural estimators is proposed in Section 4. Numerical study for the comparison of the risk values of the proposed estimators is presented in Section 5.  Section 6 concludes this paper.
	\section{{Some Inadmissibility Results\label{ch5sec5.2}}}
	The natural estimators of the hazard rate of the selected exponential population are (see 
	Jha et al. \cite{B7}) given by 
	\begin{eqnarray}
		\delta_{ML}&=&\frac{n}{Y_{J}}, \label{ch5.2.1}\\
		\delta_{N1}&=&\frac{n-2}{Y_{J}}~\mbox{and} \label{ch5.2.2}\\
		\delta_{N2}&=&\frac{n-1}{Y_{J}},
		\label{ch5.2.3}
	\end{eqnarray}
	where $\delta_{ML}$ is the maximum likelihood estimator (MLE) of $\sigma_{J}$ and $\delta_{N2}$ is the analogue of the BSEE for the component problem.
	Looking at the form of the above estimators, we consider a class of estimators of the form 
	$\delta_{c}=\frac{c}{Y_{J}},$ $c > 0,$ ($c$ is an arbitrary constant) for $\sigma_{J}$. Our aim
	is to derive a subclass of admissible estimators within this class. The risk of
	$\delta_{c}$ with respect to the loss function given by (\ref{ch5.1.1}) is
	expressed as
	$$ R(\underline{\sigma},\delta_{c}) =
	E\Bigg[\frac{\delta_{c}}{\sigma_{J}} - \ln\left(\frac{\delta_{c}}{\sigma_{J}}\right) - 1\Bigg] .$$
	It can be seen that the choice of $c$ minimizing the above risk function is found to be
	$$ c = \frac{1}{E\Big[\frac{1}{\sigma_{J} Y_{J}}\Big]}. $$
	Please note that here $c$ is a function of the parameters. To establish the following theorem, the Brewster-Zidek technique is used. For this, one needs to derive supremum and infimum of $c=\frac{1}{E\Big[\frac{1}{\sigma_{J} Y_{J}}\Big]}$. For the sake of convenience, henceforth, we denote beta function by $B(a,b) = \int_0^1 s^{a-1}(1 - s)^{b-1} ds,~  a~\textgreater~0,~  b~\textgreater~0.$ 
	\begin{thm}
		For $n ~\textgreater~ 1$ and $k = 2,$ let
		${c}^{*}=((n-1)B(n,n-1))/(2\int_{0}^{\frac{1}{2}}v^{n-2}(1-v)^{n-1}dv)$ and $c_{{*}} = n - 1$. Then, under the entropy loss function given by (\ref{ch5.1.1}), the estimator $\delta_{c}=\frac{c}{Y_{J}}$ is admissible within this class, if and only if  ${c}_{*} \leq c \leq {c}^{*}$. Further, estimators $\delta_{c}$ with $c >
			{c}^{*}$ is improved by $\delta_{c}$ with $c = {c}^{*}$ and that with $c < c_{{*}}$ is improved by $\delta_{c}$ with $c= c_{{*}}$, where ${c}^{*}$ and $c_{{*}}$ are as defined above. \label{6} 
	\end{thm}
	{\bf{Proof:}} We have
	\begin{eqnarray}
		E\Big(\frac{1}{\sigma_{J}Y_J}\Big) &=& \frac{1}{\left(\Gamma(n)\right)^2}\int_0^\infty \int_{\frac{\sigma_1 v_2}{\sigma_2}}^\infty {v_1}^{n-2}
		{v_2}^{n-1}e^{-(v_1+v_2)}dv_1 dv_2\nn\\
		&+&
		\frac{1}{\left(\Gamma(n)\right)^2}\int_0^\infty
		\int_{\frac{\sigma_2 v_1}{\sigma_1}}^\infty
		{v_2}^{n-2}{v_1} ^{n-1}e^{-(v_1+v_2)}dv_2 dv_1.
		\label{ch5neh3}
	\end{eqnarray}
The above integral can be written as
	\begin{eqnarray*}
		E\Big(\frac{1}{\sigma_{J}Y_J}\Big)&=&\frac{1}{\left(\Gamma(n)\right)^2}\int_0^\infty \int_{\frac{\sigma_1 x_2}{\sigma_2}}^\infty {x_1}^{n-2} {x_2}^{n-1} e^{-(x_1+x_2)}dx_1 dx_2\\
		&+&
		\frac{1}{\left(\Gamma(n)\right)^2}\int_0^\infty
		\int_{\frac{\sigma_2 x_2}{\sigma_1}}^\infty
		{x_1}^{n-2}{x_2} ^{n-1}e^{-(x_1+x_2)}dx_1 dx_2\\
		&=& \frac{1}{\left(\Gamma(n)\right)^2}\int_0^\infty  {x_2}^{n-1} e^{-x_2} \Bigg( \int_{\frac{\sigma_1 x_2}{\sigma_2}}^\infty {x_1}^{n-2} e^{-x_1}dx_1 \Bigg) dx_2\\
		&+&
		\frac{1}{\left(\Gamma(n)\right)^2}\int_0^\infty {x_2} ^{n-1} e^{-x_2}
		\Bigg(\int_{\frac{\sigma_2 x_2}{\sigma_1}}^\infty
		{x_1}^{n-2}e^{-x_1}dx_1 \Bigg) dx_2\\
		&=& \frac{1}{\left(\Gamma(n)\right)}\int_0^\infty  {x_2}^{n-1} e^{-x_2} \Bigg(\frac{1}{\left(\Gamma(n)\right)}  \int_{\frac{\sigma_1 x_2}{\sigma_2}}^\infty {x_1}^{n-2} e^{-x_1}dx_1 \Bigg) dx_2\\
		&+&
		\frac{1}{\left(\Gamma(n)\right)}\int_0^\infty {x_2} ^{n-1} e^{-x_2}
		\Bigg( \frac{1}{\left(\Gamma(n)\right)} \int_{\frac{\sigma_2 x_2}{\sigma_1}}^\infty
		{x_1}^{n-2}e^{-x_1}dx_1 \Bigg) dx_2.\\
	\end{eqnarray*}
	Hence, $E\Big(\frac{1}{\sigma_{J}Y_J}\Big)$
	\begin{eqnarray*}
		&=& \frac{1}{\left(\Gamma(n)\right)}\int_0^\infty  {x_2}^{n-1} e^{-x_2} \Bigg(\frac{1}{(n - 1)\left(\Gamma(n - 1)\right)}  \int_{\frac{\sigma_1 x_2}{\sigma_2}}^\infty {x_1}^{n-2} e^{-x_1}dx_1 \Bigg) dx_2\\
		&+&
		\frac{1}{\left(\Gamma(n)\right)}\int_0^\infty {x_2} ^{n-1} e^{-x_2}
		\Bigg( \frac{1}{(n - 1)\left(\Gamma(n - 1)\right)} \int_{\frac{\sigma_2 x_2}{\sigma_1}}^\infty
		{x_1}^{n-2}e^{-x_1}dx_1 \Bigg) dx_2\\
		&=& \frac{1}{\left(\Gamma(n)\right)} \frac{P(X_1>q X_2)}{(n - 1)} \int_0^\infty  {x_2}^{n-1} e^{-x_2}  dx_2\\
		&+&
		\frac{1}{\left(\Gamma(n)\right)} \frac{P\Big(X_1>\frac{X_2}{q}\Big)}{(n - 1)} \int_0^\infty {x_2} ^{n-1} e^{-x_2} dx_2\\
		&=& \frac{1}{n-1}\Big[P(X_1>q X_2)+P(X_1>\frac{X_2}{q})\Big].\\  
\end{eqnarray*}

\noindent Suppose we have two independent gamma random variables $X_1$ and $X_2$ having a common scale parameter $1$ and different shape parameters $n-1$ and $n$, respectively. Let $V=\frac {X_2}{X_1+X_2}.$
Then, $V$ follows beta distribution with parameters $n$ and $n-1$. Define $q = \frac{\max(\sigma_1,\sigma_2)}{\min(\sigma_1,\sigma_2)} ( \geq 1).$ Again, (\ref{ch5neh3}) can be written as
\begin{eqnarray}
		E\Big(\frac{1}{\sigma_{J}Y_J}\Big)&=& \frac{1}{n-1}\Bigg[P(X_1>q X_2)+P(X_1>\frac{X_2}{q})\Bigg]\nn\\
		&=& \frac{1}{n-1}\Bigg[P\Big(\frac{X_1}{X_2} + 1 > q + 1\Big)+ P\Big(\frac{X_1}{X_2} + 1>\frac{1}{q} + 1\Big)\Bigg]\nn\\
		&=& \frac{1}{n-1}\Bigg[P\Big(\frac{X_2}{X_2+X_1}<\frac{q}{1+q}\Big)+ P\Big(\frac{X_2}{X_2+X_1}<\frac{1}{1+q}\Big)\Bigg]\nn\\
		&=& \frac{\Bigg[\int_0^{\frac{q}{1+q}}v^{n-1}(1-v)^{n-2}dv+
			\int_0^{\frac{1}{1+q}}v^{n-1}(1-v)^{n-2}dv\Bigg]}{(n-1)B(n,
			n-1)}\nn\\
		& =& h(q)~(say).\label{ch5neh4}
\end{eqnarray} 
	Differentiating $h(q)$ with respect to $q$, we obtain
	\begin{eqnarray*}
		h'(q)& =&  \frac{q^{n-2}}{(n-1)B(n,n-1)}\Big[\frac{(q-1)}{(1 + q)^{2n-1}}\Big] \\
		&\geq& 0, 
	\end{eqnarray*}
	where the inequality holds for $q \geq 1$.
	Thus, $h(q)$ is increasing for $q \in [1,\infty)$. This implies that $c$ is decreasing with respect to $q$ in  $[1,\infty).$ Now, we obtain the supremum and infimum of $h(q)$ as follows:
	
	\begin{eqnarray*}
		\sup_q c =  {c}^*
		=\lim_{q\rightarrow 1}c
		=\lim_{q\rightarrow 1}\frac{1}{h(q)}
		=\frac{(n-1)B(n,n - 1)}
		{2\int_0^{\frac{1}{2}}v^{n-1}(1-v)^{n-2}dv}\\
	\end{eqnarray*}
	and
	\begin{eqnarray*}
		\inf_q c= {c}_{*}=\lim_{q\rightarrow
			\infty}c &=&\lim_{q\rightarrow \infty}{\frac{1}{h(q)}}\\
		&=& n - 1. \label{neh8}
	\end{eqnarray*}
	Applying Brewster-Zidek technique (see Theorem $3.3.1$, page-$36$ of Brewster and Zidek \cite{B5} for details in this regard), the estimator $\delta_{c}$ is admissible within the class of estimators of the form $\delta_{c}=\frac{c}{Y_{J}}$, where ${c}_{*} \leq c \leq {c}^{*}$ is the range of admissibility. Hence, the proof is completed.\\
	
	The following remark is an immediate consequence of Theorem $2.1$.
	\begin{rem}
		Since $(n - 2)$~\textless~$(n - 1),$ the estimator $\delta_{N1}= \frac{n-2}{Y_{J}}$ is inadmissible. And, the estimator $\delta_{N2}$ dominates $\delta_{N1}$ according to Theorem \ref{6}.
	\end{rem}

	\section{{A Minimaxity Result\label{ch5sec5.3}}}
	In this section, we establish minimaxity of $\delta_{N2}$. To do this, let us consider the estimation of the hazard rate $\sigma$ for the component problem. A component problem is the one in which selection criteria is not taken care of. Bayes and generalized Bayes estimators are obtained in Lemma $\ref{3.1}$. The proof of the Bayes and generalised Bayes estimators is quite straightforward, and hence it is omitted. However, the proof for other part of the Lemmas have been provided.
	\begin{lem}
		Let $Y$ follow Gamma($\sigma,n$) distribution. Then, under the entropy loss function,
		\begin{itemize}
			\item[I.] the Bayes estimator of $\sigma$ with respect to the conjugate gamma prior
			$\pi\sim Gamma(\alpha,\gamma)$ is given by $\delta_{B}=\frac{n+\alpha-1}{Y+\gamma}$.
			The Bayes risk of $\delta_{B}$ is $r(\sigma,\delta_{B})=\Psi(n + \alpha)- \ln(n + \alpha - 1)$, where $\Psi(.)$ is a digamma function;
			\item[II.] a generalized Bayes estimator with respect to the noninformative prior
			$\tau(\sigma)=1/\sigma,~\sigma>0$ is obtained as $\delta_{GB}=\frac{n-1}{Y}$. The risk of $\delta_{GB}$ is $R(\sigma,\delta_{GB})=\Psi(n)- \ln(n - 1)$. Further, the estimator $\delta_{GB}=\frac{n-1}{Y}$ is minimax and the minimax value is $\Psi(n)- \ln(n - 1)$.\\
			
			\noindent {\bf{Proof:}} 
			
			I. The Bayes risk of $\delta_B$ with respect to the entropy loss function is given by
			\begin{eqnarray*}
				r(\sigma,\delta_B)  &=& EE\Bigg[\frac{n + \alpha -1}{\sigma(Y + \gamma)} - \ln\frac{n + \alpha -1}{\sigma(Y + \gamma)} - 1\Bigg]\\
				&=& \int_{0}^{\infty} \int_{0}^{\infty} \Bigg[\frac{n + \alpha -1}{\sigma(y + \gamma)} - \ln\frac{n + \alpha -1}{\sigma(y + \gamma)} - 1\Bigg]f(\sigma,y)d\sigma dy\\
				&=&  \int_{0}^{\infty} \int_{0}^{\infty} \Bigg[\frac{n + \alpha -1}{\sigma(y + \gamma)} - \ln\frac{n + \alpha -1}{\sigma(y + \gamma)} - 1\Bigg]\frac{\gamma^{\alpha}}{\Gamma(\alpha)\Gamma(n)}~ y^{n-1}~\sigma^{n+\alpha-1}~ e^{-\sigma(y+\gamma)}dyd\sigma.
			\end{eqnarray*}
			Proceeding further, finally we get
			\begin{eqnarray*}
				r(\sigma,{\delta_B})
				&=& \frac{\gamma^{\alpha}}{\Gamma(\alpha)\Gamma(n)}\Bigg[\Gamma^{'}(n + \alpha) - \Gamma(n + \alpha)\ln(n + \alpha -1) \Bigg] \int_{0}^{\infty} \frac{y^{n - 1}}{(y + \gamma)^{n + \alpha}}dy\\
				&=& \Psi(n + \alpha) - \ln(n + \alpha -1),
			\end{eqnarray*}
		    where $\Gamma^{'}(n + \alpha) = \int_{0}^{\infty} \ln t ~ t^{n + \alpha -1} ~ e^{-t} dt$.\\
		    
			 II. We know that,
			\begin{eqnarray*}
				R(\sigma,\delta_{GB}) &=& R\Big(\sigma,\frac{n - 1}{Y}\Big)\\
				 &=&  E \Big[\Big(\frac{n - 1}{\sigma Y}\Big) - \ln\Big(\frac{n - 1}{\sigma Y}\Big) - 1\Big]\\
				&=&  (n - 1) E\Big(\frac{1}{\sigma Y}\Big) - \ln(n - 1) + E(\ln(\sigma Y)) - 1,
			\end{eqnarray*}
			where~
			\begin{eqnarray*}
				E\Big(\frac{1}{\sigma Y}\Big)  &=& \int_{0}^{\infty} \frac{1}{\sigma y} \frac{\sigma^{n}}{\Gamma(n)}  ~ y^{n-1} ~e^{-\sigma y} dy\\
				&=&  \frac{\sigma^{n - 1}}{\Gamma(n)} \int_{0}^{\infty} y^{n-2} ~e^{-\sigma y} dy\\
				&=&  \frac{\Gamma(n - 1)}{\Gamma(n)}\\
				&=& \frac{1}{n - 1}
			\end{eqnarray*}
			and~
			\begin{eqnarray*}
				E(\ln(\sigma Y))  &=&  \int_{0}^{\infty} \ln(\sigma Y)  \frac{\sigma^{n}}{\Gamma(n)}  ~ y^{n-1} ~e^{-\sigma y} dy\\
				&=&  \frac{\sigma^{n}}{\Gamma(n)\sigma^{n}} \int_{0}^{\infty} \ln(t)  ~ t^{n-1} ~e^{-t} dt\\
				&=&  \frac{\Gamma^{'}(n)}{\Gamma(n)}\\
				&=&  \Psi(n). 
			\end{eqnarray*}
			Hence,
			\begin{eqnarray*}
				R\Big(\sigma,\frac{n - 1}{Y}\Big) &=& (n - 1)\Big(\frac{1}{n - 1}\Big) - \ln(n - 1) + E(\ln(\sigma Y)) - 1\\
				&=& 1 - 1 + \Psi(n) - \ln(n - 1)\\
				&=&\Psi(n) - \ln(n - 1).
			\end{eqnarray*}
		\end{itemize}
		\label{3.1} 
		Hence, the risk value of the estimator $\delta_{GB} = \frac{n - 1 }{Y}$ equals $\Psi(n) - \ln(n - 1)$. Further, $R(\sigma,\delta_{GB})$ = $\lim\limits_{\alpha \rightarrow 0, \gamma \rightarrow 0}r(\sigma,\delta_B)$. So, the estimator $\delta_{GB}$ is minimax.
	\end{lem}
	
	Next, we prove a theorem, which deals with minimaxity of $\delta_{N2}$.

	\begin{thm}
		The estimator $\delta_{N2} = \frac{n - 1}{Y_{J}}$ is minimax for the estimation of $\sigma_{J}$ with respect to the loss function given by (\ref{ch5.1.1}).
	\end{thm}
	{\bf{Proof:}}~ 
	Following Sackrowitz and Samuel Cahn \cite{B15}, we proceed as follows
	\begin{eqnarray*}
		R\Bigg(\underline{\sigma},\frac{n - 1}{Y_{J}}\Bigg) &=& E \Bigg[\Bigg(\frac{n - 1}{\sigma_{J}Y_J}\Bigg) - \ln\Bigg(\frac{n - 1}{\sigma_{J}Y_J}\Bigg) - 1\Bigg]\\
		&=& (n - 1)E\Bigg(\frac{1}{\sigma_{J}Y_J}\Bigg) - \ln(n - 1) + E(\ln (\sigma_{J}Y_{J})) - 1\\
	\end{eqnarray*}
	Proceeding similarly as in Theorem $2.1$, we have\\
		$$\sup E(\ln(\sigma_{J}Y_J)) \leq \Psi(n) - \Psi(2n - 1) \leq \Psi(n) \quad \forall ~n \geq 2.$$
	Hence,
	\begin{eqnarray*}
		R\Bigg(\underline{\sigma},\frac{n - 1}{Y_{J}}\Bigg)   	&\leq& (n - 1)h(q) - \ln(n - 1) + \Psi(n) - 1
		= H(q)~ \mbox{(say)}, \mbox{where}~ q\geq1.
	\end{eqnarray*}
	
	\noindent So, $$H^{'}(q) = (n - 1)h^{'}(q),$$
	which is non-negative since $n> 1$ and  $h^{'}(q) \geq 0.$
	Hence, $H(q)$ is increasing. Further,  
	\begin{eqnarray*}
		\sup_q H(q) = \lim\limits_{q \rightarrow \infty} H(q)
		=   \Psi(n) -\ln(n - 1).
	\end{eqnarray*}
	This implies that $H(q) \leq \Psi(n) - \ln(n - 1)$.
	According to Theorem $3.2$ of Sackrowitz and Samuel-Cahn \cite{B15}, an estimator $\delta$ of $\sigma_{J}$
	is minimax, if
	$\sup_{\underline{\sigma}}
	R(\underline{\sigma},\delta ) \leq \sum_{i=1}^{k}
	r_i^{*}(\delta_B)I_i$, where
	$r_i^{*}(\delta_B,x_i)=r_i^{*}(\delta_B)$ is the
	posterior risk of the Bayes estimator $\delta_B$
	of $\sigma_i$ in the $i$th component estimation
	problem.
	Therefore, we conclude that $R(\underline{\sigma}, \delta_{GB}) \leq \Psi(n) - \ln(n - 1).$ This proves the theorem.
	
	%
	%
	
	\begin{rem}
		Proceeding similar to above, we see that
		$\sup_{\underline{\sigma}}
		R(\underline{\sigma},\delta_{N1}) = \Psi(n) - \ln(n - 2) - \frac{1}{n - 1}$ and $\sup_{\underline{\sigma}}
		R(\underline{\sigma},\delta_{ML}) =  \Psi(n) - \ln(n) + \frac{1}{n - 1},$ which are greater than the
		minimax value $\Psi(n) - \ln(n - 1).$ Hence, the estimators
		$\delta_{N1}$ and $\delta_{ML}$ are not minimax.
	\end{rem}
	
	In the next section, an emphasis is
	made to obtain some estimators which improve upon
	the natural estimators in terms of their risk
	values.
	
	\section{{Improved Estimators\label{ch5sec5.4}}}
	
	We would like to obtain new estimators
	which dominate the  natural estimators as 
	proposed in Section \ref{ch5sec5.2}. First, we quote a lemma (see Berger \cite{B4} and Vellaisamy \cite{B18}), which is useful to prove some of the results in this section. Let $\underline{\alpha} = (\alpha_1,\ldots,\alpha_k).$
	
	\begin{lem} (Vellaisamy \cite{B18})
		Let $X_1,\ldots,X_k$ be $k$ independent random variables, where $X_i$ follows $Gamma(\alpha_i, p_i)$ distribution, $i= 1,\ldots,k.$
		Let $U(\underline{x})$ be any real valued function on $\mathbb{R}^{k}$ such that\\ $(i)~E_{\underline{\alpha}} \left|U(\underline{X})\right|<\infty$, where $\underline{X} = (X_1,\ldots,X_k)$ and $\underline{x} = (x_1,\ldots,x_k)$
		for all $\underline{\alpha}$ and\\ (ii) the definite integral
		\begin{eqnarray*}
			h_i(\underline{x}) = \int_{0}^{x_i} U(x_1,\ldots,x_{i-1},t,x_{i+1},\ldots,x_k)t^{p_i-1}dt
		\end{eqnarray*}
		\quad exists for all $x_i \in \mathbb{R}$. Then, $V(\underline{x}) = x_i^{1 - p_i}h_i(\underline{x})$ satisfies the condition
		$$E_{\underline{\alpha}}(V(\underline{X})) = \alpha_i E_{\underline{\alpha}}(U(\underline{X})),~  for~ all~ \underline{\alpha}.$$ \label{lem1}
	\end{lem}
	
	It may be noted that the above lemma is used to derive the unbiased estimator of the risk difference. For the purpose of deriving new estimators, let $\underline{Y}=(Y_1,\ldots,Y_k).$ Suppose  $Y_{(1)i} \geq Y_{(2)i} \geq \ldots \geq Y_{(k - 1)i}$ denote the ordered values of $(Y_1,\ldots,Y_{i - 1}, Y_{i + 1},\ldots,Y_k)$. 	Next, we prove a theorem to find out the improved estimators over the natural estimators for estimating $\sigma_{J}$ which is the hazard rate of the selected population.
	
	

	\begin{thm}
		Let $X = \Big(\prod_{i = 1}^{k}Y_{i}\Big)^{\frac{1}{k}}$ and $w(t), ~t>  0$ be a real valued function such that
		\begin{itemize}
			\item[(i)] $ w(t)~is~ nonincreasing$ in $(0,\infty)$ and
			\item[(ii)]  $0 ~\textless~tw(t) \leq \frac{(n - c)k + 1}{nk + 1}$, $t~\textgreater~0.$
		\end{itemize}
		Then, for each $0~\textless~c\leq n,$ any estimator of the form\\
			$$\delta_{ri}= 
			\Bigg[\frac{c}{Y_{i}} + n\Big[w(X) + \frac{Xw'(X)}{nk}\Big]\Bigg],$$
			dominates the estimator $\frac{c}{Y_{i}}.$ \label{4.1}\end{thm}
	{\bf{Proof:}}~ 
	Let $\delta_2(\underline{Y})  = \Big(\delta_{21}I_1,\ldots,\delta_{2k}I_k\Big),$ be any improved estimator of $\underline{\sigma}$ improving upon \\ $\delta_1(\underline{Y})  =  \Big(\delta_{11}I_1,\ldots,\delta_{1k}I_k\Big).$
	The risk difference of $\delta_1(\underline{Y})$ and $\delta_2(\underline{Y})$ for estimating $\underline{\sigma}$ under the loss function (\ref{ch5.1.1}) is obtained as
	\begin{eqnarray*}
		R(\delta_2(\underline{Y}),\underline{\sigma}) - R(\delta_1(\underline{Y}),\underline{\sigma})
		&=&E\sum_{i=1}^{k}\Bigg[\Bigg(\frac{\delta_{2i}}{\sigma_i} - \ln\frac{\delta_{2i}}{\sigma_i} -1\Bigg) - \Bigg(\frac{\delta_{1i}}{\sigma_i} - \ln\frac{\delta_{1i}}{\sigma_i} -1\Bigg)\Bigg] I_i\nn\\
		&=& E\sum_{i =1}^{k}\Bigg[\frac{1}{\sigma_i}{(\delta_{2i} - \delta_{1i})} +  \ln\Bigg(\frac{\delta_{1i}}{\delta_{2i}}\Bigg)\Bigg]I_i.\nn\\
	\end{eqnarray*}
	After applying Lemma \ref{lem1} to the first term inside the square brackets in the above expression, the unbiased estimator of the  risk difference is derived as
	\begin{equation}
		R(\delta_2(\underline{Y}),\underline{\sigma}) - R(\delta_1(\underline{Y}),\underline{\sigma})  = E\sum_{i=1}^{k}\Bigg[ Y_i^{1 - n}\int_{Y_{(1)i}}^{Y_i}(\delta_{2i}(\underline{Y})-\delta_{1i}(\underline{Y}))t^{n - 1}dt + \ln\Bigg(\frac{\delta_{1i}(\underline{Y})}{\delta_{2i}(\underline{Y})}\Bigg) \Bigg]I_i. \label{1}\end{equation}
	Let $G_i(\underline{Y}) = \frac{1}{n} \int_{Y_{(1)i}}^{Y_i}(\delta_{2i}(\underline{Y})-\delta_{1i}(\underline{Y}))t^{n - 1}dt$ such that $G_i^{i(1)}(\underline{Y})= \frac{1}{n} (\delta_{2i}(\underline{Y})-\delta_{1i}(\underline{Y}))y_i^{n - 1},$ where $G_i^{i(1)}(\underline{Y})$ is the first order partial derivative of $G_i(\underline{Y})$ with respect to $y_i.$
	Substituting these values in (\ref{1}), we get 
	$$ R(\delta_2(\underline{Y}),\underline{\sigma}) - R(\delta_1(\underline{Y}),\underline{\sigma})= E\sum_{i=1}^{k}\Bigg[Y_{i}^{1 - n}nG_i(\underline{Y}) + \ln \Bigg(\frac{\delta_{1i}(\underline{Y})}{\delta_{2i}(\underline{Y})}\Bigg)\Bigg]I_i.$$
	Suppose, we take $\delta_{1i}=\frac{c}{Y_i}$ and $\delta_{2i}=\delta_{1i} + Y_{i}^{1 - n}nG^{i(1)}(\underline{Y})$.
	Now, (\ref{1}) can be expressed as
	\begin{eqnarray}
		R(\delta_2(\underline{Y})) - R(\delta_1(\underline{Y}))
		= E\sum_{i=1}^{k}\Bigg[Y_{i}^{1 - n}(nG(Y_{1},\ldots,Y_{k})) - \ln \Big(1 + \frac{n}{c}Y_{i}^{2 - n} G^{i(1)}(Y_{1},\ldots,Y_{k})\Big)\Bigg]I_i.\nn\\
		\label{7} \end{eqnarray}
	We have to find out a suitable function $G(y_{1},\ldots,y_{k})$ such that (\ref{7}) is negative in the region $B = \{\underline{y}:y_{1}\geq y_{2}\geq\ldots\geq y_{k} ~ \textgreater ~ 0\}$, i.e, we aim to solve the differential inequality $Q(y_{1},\ldots,y_{k}) \leq 0$ in the above mentioned region, where\\
	\begin{equation}
		Q(y_{1},\ldots,y_{k}) = \Bigg[y_{1}^{1 - n}(nG(y_{1},\ldots,y_{k})) - \ln \Big(1 + \frac{n}{c} y_{1}^{2 - n} G^{i(1)}(y_{1},\ldots,y_{k})\Big)\Bigg]. \label{8}
	\end{equation}
	We see that whenever $\underline{y} \in B$, then clearly $y_{1} \geq x,$ when $x = (\prod_{i = 1}^{k}y_{i})^{\frac{1}{k}}.$
	Now, let us choose\\
	$$G(y_{1},\ldots,y_{k}) = \frac{y_{1}^n w(x)}{n}, ~\mbox{for some}~ w(x),$$\\
	where $w(x)$ satisfies the assumptions of the theorem. It is observed that
	\begin{eqnarray*}
		G^{i(1)}(y_{1},\ldots,y_{k}) = y_{1}^{n - 1}\Big[w(x) + \frac{xw'(x)}{nk}\Big]~ \mbox{and} ~\Big[w(x) + \frac{xw'(x)}{nk}\Big] \leq \Big(1 + \frac{1}{nk}\Big)w(x).
	\end{eqnarray*}
	Now, substituting the values of $G_i(\underline{y})$ and $\delta_{1i}$ in the above expression for $\delta_{2i}$ above, we get
	\begin{eqnarray*}
		\delta_{2i}(y_{1},\ldots,y_{k}) = \frac{c}{{y_{1}}
		} + n\Big[w(x) + \frac{xw'(x)}{nk}\Big].\label{4}
	\end{eqnarray*} 
	Hence (\ref{8}) becomes,
	\begin{eqnarray*}
		Q(y_{1},\ldots,y_{k}) &=&\Bigg[y_{1} w(x) - \ln\Bigg(1 + \frac{n}{c}y_{1} \Big(w(x) + \frac{xw'(x)}{nk}\Big)\Bigg) \Bigg],\\
		&\leq&  \Bigg[y_{1} w(x) - \frac{\frac{n}{c}{y_{1}} \Big(w(x) + \frac{xw'(x)}{nk}\Big)}{1 + \frac{n}{c}y_{1} \Big(w(x) + \frac{xw'(x)}{nk}\Big)} \Bigg] \\
		&\leq& \Bigg[\frac{y_{1} w(x)\Big(1 + \frac{n}{c}y_{1} \Big(w(x) + \frac{xw'(x)}{nk}\Big)\Big) - \frac{n}{c}y_{1} \Big(w(x) + \frac{xw'(x)}{nk}\Big)} {1 + \frac{n}{c}y_{1} \Big(w(x) + \frac{xw'(x)}{nk}\Big)} \Bigg]\\
		&\leq& \Bigg[y_{1} w(x) \Big(\frac{c - n}{c} \Big)+ y_{1}^2w(x)^2 \Big(\frac{kn + 1}{ck}\Big) - \frac{{y_{1}} xw'(x)}{ck}\Bigg]\\
		&\leq& 0,
	\end{eqnarray*}
	where the first inequality is due to the fact that $\ln(1 + y) \geq \frac{y}{1 + y}$,  for $ y>-1$.
	The above inequality is satisfied, if
	\begin{eqnarray*}
		&~&y_{1} w(x) \Big(\frac{c - n}{c} \Big)+ y_{1}^2w(x)^2 \Big(\frac{kn + 1}{ck}\Big)\leq ~\frac{{y_{1}} xw'(x)}{ck}\leq \frac{{y_{1}} w(x)}{ck}\\
		&~\Rightarrow ~& xw(x) \leq \frac{(n - c)k + 1}{kn + 1}.
	\end{eqnarray*}  Hence, it can be shown that $G$ chosen  above is a solution of the differential inequality $Q(y_{1},\ldots,y_{k}) \leq 0$. Thus, the theorem is proved.\\

	Next, one would like to obtain $w(t)$ which satisfies all the conditions of Theorem \ref{4.1}.
	\begin{rem}
		Let $w(t)= \frac{\alpha}{t}, ~\mbox{where} ~ 0~\textless~\alpha \leq \frac{(n - c)k + 1}{nk + 1},~ for~ 0~\textless~c\leq n.$ This choice of $w(t)$ satisfies all the conditions of Theorem \ref{4.1}. \label{rem 1}
	\end{rem}
	\begin{rem}
		We can also take X as the geometric mean of the first $h$ order statistics, that is $X = (\prod_{i = 1}^{h}Y_{(i)})^{\frac{1}{h}}, 2 \leq h \leq k.$ This theorem still holds good for the upper bound $\frac{((n - c)h + 1)}{(nh + 1)}.$ Applying Theorem \ref{4.1}, we can derive the improved estimators of $\delta_{N2}$ and $\delta_{ML}$. It may be noted that, the estimator $\frac{n - 2}{Y_J}$ is improved by $\frac{n - 1}{Y_J}$ using Brewster-Zidek technique as discussed in Section \ref{ch5sec5.2}. \label{rem 2}\end{rem}
	
	Further,  some results on improvement of natural estimators except $\delta_{N1}$(as discussed in $Remark~ \ref{rem 2}$) is straightforward from $Theorem~ \ref{4.1}$ and $Remark~\ref{rem 1},$ as stated in the following manner.
	\begin{cor}
		For a suitably chosen $w(x)$, if $c = n - 1,$ then upper bound of $xw(x)~\mbox{is}~\frac{k + 1}{kn + 1}.$
		For $n > 1$, putting $c = n - 1,$ the estimator $\delta_{N2}=\frac{n - 1}{Y_{J}}$ is dominated by $\delta_{N2}^{I}=\frac{n - 1}{Y_{(1)}} + n\Big[w(X) + \frac{Xw'(X)}{nk}\Big],$ with respect to the loss function (\ref{ch5.1.1}).
	\end{cor}
	\begin{rem}
		Choosing $w(x)= \frac{\alpha}{x},$ the estimator improving upon $\delta_{N2}$ is given by $\delta_{N2}^{I}=\frac{n - 1}{Y_{(1)}} + \Big[\frac{\alpha(nk - 1)}{Xk}\Big], for ~ 0 ~\textless~ \alpha \leq \frac{k + 1}{kn + 1}.$
	\end{rem}
	\begin{cor}
		For a suitably chosen $w(x)$, if $c = n,$ the upper bound of $xw(x)~\mbox{is}~\frac{1}{kn + 1}.$ Putting $c = n,$ the estimator $\delta_{ML}=\frac{n}{Y_{J}}$ is dominated by $\delta_{ML}^{I}=\frac{n}{Y_{(1)}} + n\Big[w(X) + \frac{Xw'(X)}{nk}\Big],$ with respect to the loss function (\ref{ch5.1.1}).
	\end{cor}
	\begin{rem}
		Choosing $w(x)= \frac{\alpha}{x},$ the estimator improving upon $\delta_{ML}$ is given by  $\delta_{ML}^{I}=\frac{n}{Y_{(1)}} + \Big[\frac{\alpha(nk - 1)}{Xk}\Big], for ~ 0 ~\textless~ \alpha \leq \frac{1}{kn + 1}.$
	\end{rem}
	In the next section, we numerically compare the risk values of various established estimators.
	\section{Numerical Comparisons}\label{ch5sec5.5}
	In this section, simulation is done to have a comparative study of the proposed estimators in terms of risk values. For this purpose, five thousand random samples are generated from the one-parameter exponential distribution for different choices of the sample sizes and scale parameter $\frac{1}{\sigma_i}.$ The simulation study has been performed using $R$ software for $k = 2$ populations. The risk values of the natural estimators $\delta_{N2}$, $\delta_{N1}$ and  $\delta_{ML}$ with respect to the entropy loss function have been computed, tabulated and compared in the following tables. We have also calculated the risk values for the improved estimators $\delta_{N2}^I$ and $\delta_{ML}^I$. A homogeneous trend is observed for other sample sizes and hence omitted. We have presented the tabulated risk values graphically only for estimators $\delta_{N2}$ and $\delta_{ML}$ since $\delta_{N2}$ dominates $\delta_{N1}$ in terms of risk values, so the risk analysis of the estimator $\delta_{N1}$ is omitted. For the sake of convenience, we have plotted the values of $\Big(\frac{\sigma_1^{-1}}{\sigma_2^{-1}}\Big)$ along the x-axis and risk values along the y-axis in the graphs. From the simulation study, the following points have been observed.
	\begin{itemize}
		\item[(i)] The risks of the estimators decrease when $n$ increases. 
		\item[(ii)] From the simulation results, the estimator $\delta_{N2}$ always dominates estimator
		$\delta_{N1}$, which is in accordance with the results obtained from Brewster-Zidek technique in Section \ref{ch5sec5.2}. However, we are not able to obtain a  hierarchy between the estimators $\delta_{N2}$ and $\delta_{ML}$ as they are not comparable. 
		\item[(iii)] The margin of the improvement decreases as the sample size increases.
		\item[(iv)] It is observed that $\delta_{N2}$ has a better risk improvement over a major part of the parameter space.
	\end{itemize}
	\newpage
	\small{
		\hspace{1.3cm}{\bf Table 5.1: $n=5~~~~~~~~~~~~~~~~~~~~~~~~~~~~~$}
		\begin{center}
			\begin{tabular}{|c|c||c|c|c||c|c||}
				\hline
				\emph{$\sigma_1^{-1}
					\downarrow$}&\emph{$\sigma_2^{-1}
					\downarrow$}&{\emph{$R(\delta_{N1})$}}&{\emph{$R(\delta_{N2})$}}&{\emph{$R(\delta_{N2}^I)$}}
				&{\emph{$R(\delta_{ML})$}}&{\emph{$R(\delta_{ML}^I)$}}\\
				\hline
				~&0.2&	0.20076&  0.106062&  0.101426&  0.075902&  0.074986\\
				~&0.4&	0.206475&  0.106257&	0.101412&  0.070577& 0.069361\\
				0.3&0.6&	0.179391&  0.098671&	0.094465&  0.082489& 0.082402\\
				~&0.8&	0.165354&  0.103346&	0.099937&  0.105877& 0.104411\\
				~&1&	0.155402&  0.100987&  0.09768&	  0.111111&	0.098443\\
				\hline
				~&0.2&	0.171395&  0.10129&	0.097336&  0.095723&	0.094488\\
				~&0.4&	0.211737&  0.109796&	0.104847&  0.072393&0.071058\\
				0.5&0.6&	0.209921&  0.106511&	0.101504&  0.067639& 0.066228\\
				~&0.8&	0.199194&  0.104464&	0.099708&  0.074272& 0.073331\\
				~&1&	0.18223&  0.100169&	0.095909&  0.082647&	0.077492\\
				\hline
				\hline
				~&0.2&	0.159905&  0.103137&	0.099437&  0.110908&	0.109972\\
				~&0.4&	0.18869&	  0.100698&	0.096247&  0.077244&	0.076729\\
				0.7&0.6&	0.214956&  0.110352&	0.105241&  0.070287&	0.068772\\
				~&0.8&	0.215049&  0.110853&	0.105805&  0.071195&	0.069729\\
				~&1&	0.201175&  0.103381&	0.098612& 0.070125&	0.066812\\
				\hline
				~&0.2&	0.153211&  0.108142&	0.104994&  0.127612&  0.111144\\
				~&0.4&	0.173765& 0.100596&	0.096725&  0.091966&	0.090528\\
				0.9&0.6&	0.200779& 0.104491&	0.099718&		0.072742&	0.071736\\
				~&0.8&	0.214185& 0.109479&	0.104392&	0.069312& 0.067798\\
				~&1&	0.211906& 0.107453&	0.102429&		0.067539&	0.063751\\
				\hline
				~&0.2&	0.155423& 0.110459&	0.107136&		0.130032&	0.129764\\
				~&0.4&	 0.168475& 0.099424&	0.095599&		0.094912& 0.093843\\
				1&0.6&	0.19756&	 0.106854&	0.102241&		0.080687&	0.079955\\
				~&0.8&		0.211814& 0.108805&	0.103746&	0.070334&	0.068904\\
				~&1&	0.215238& 0.10847&	0.103268&		0.066240&	0.061222\\
				\hline
				\hline
			\end{tabular}
		\end{center}
	}
	\begin{figure}[h!]
		\begin{subfigure}{0.6\textwidth}
			\centering
			\includegraphics[width=0.7\linewidth]{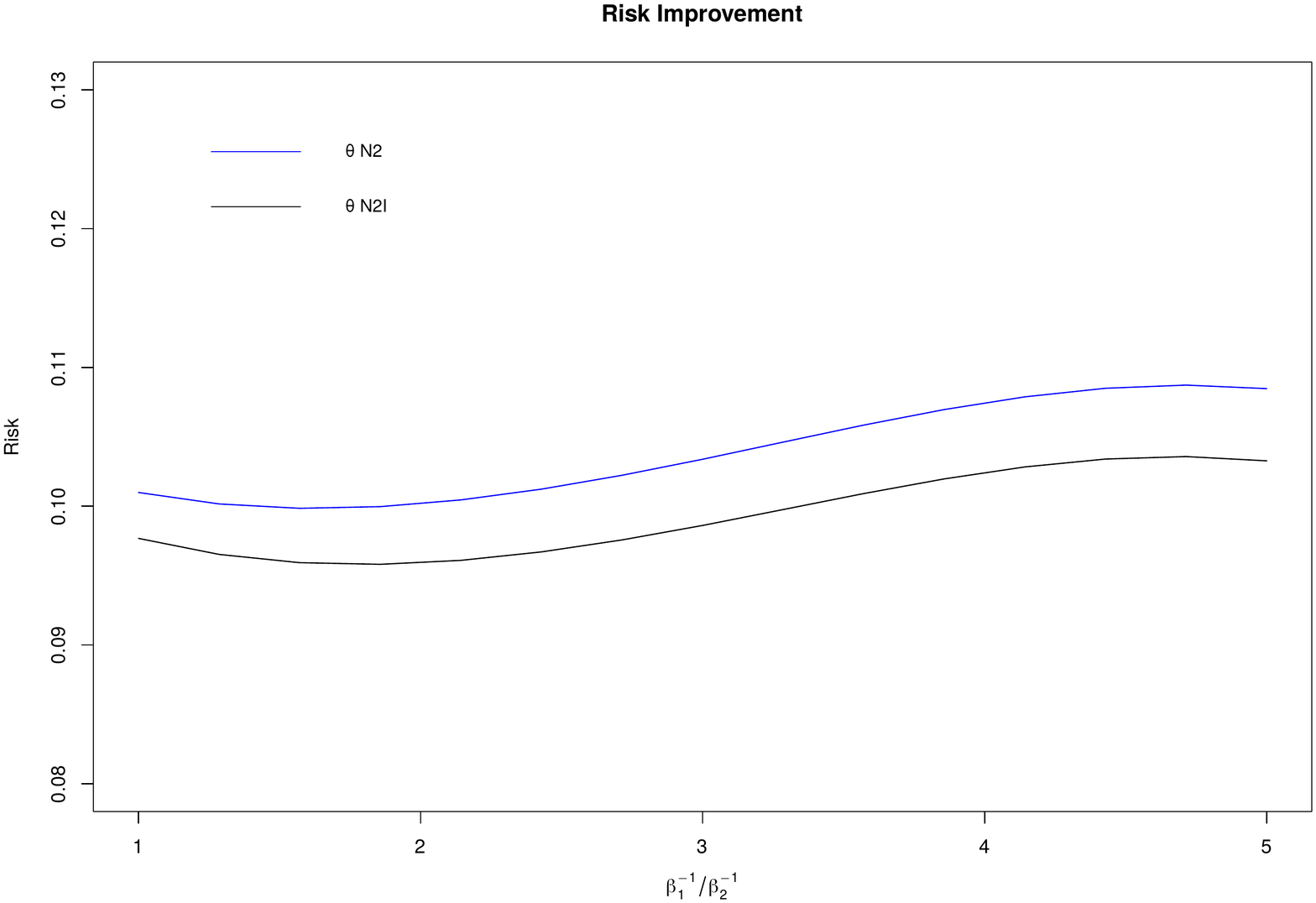}
			\caption{}
			\label{fig:test1}
		\end{subfigure}%
		\begin{subfigure}{0.6\textwidth}
			\includegraphics[width=0.7\linewidth]{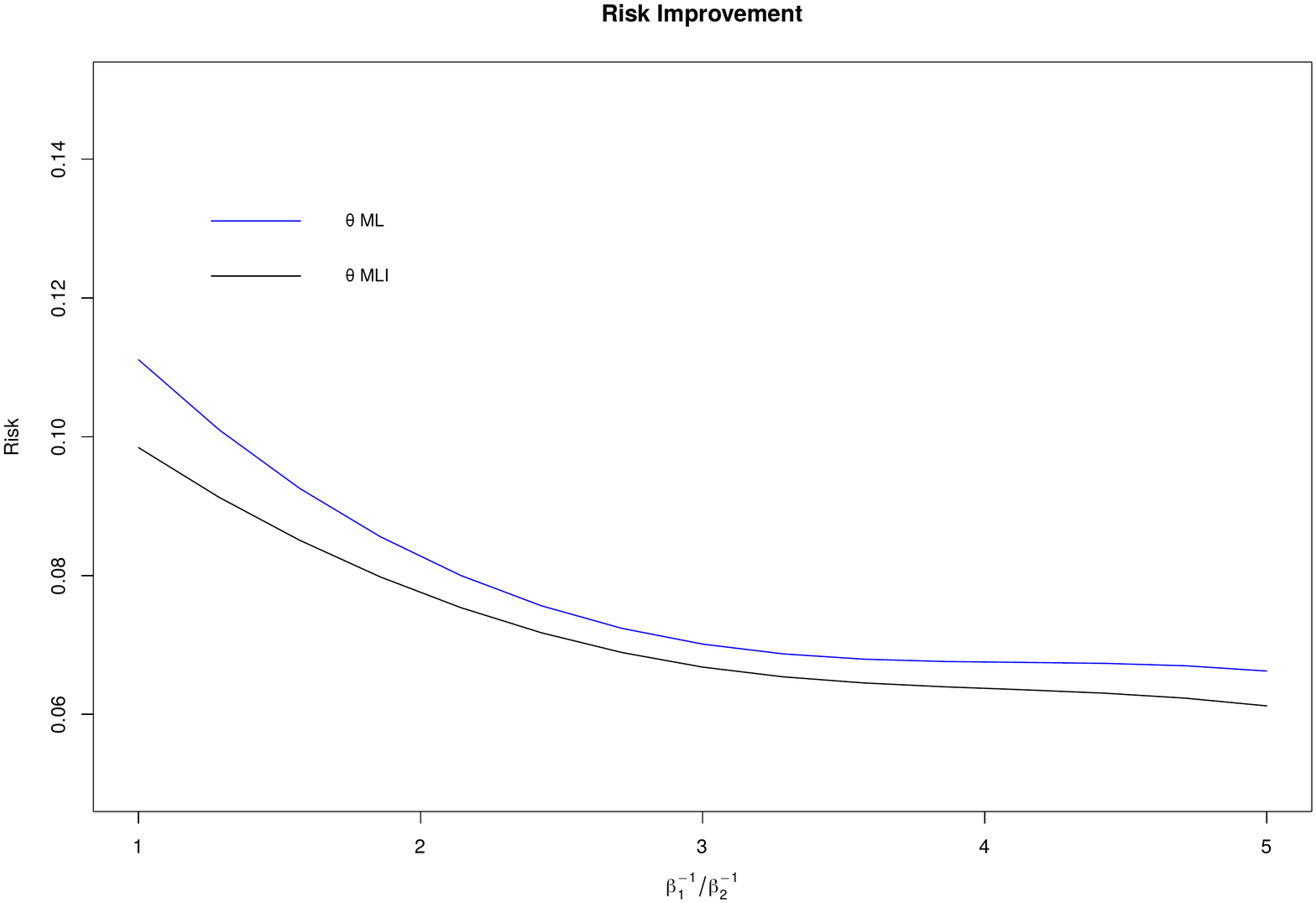}
			\caption{}
			\label{fig:test2}
		\end{subfigure}
		\caption{ For $n=5,$ the risk plots of the estimators 
			(a)$\delta_{N2}$, $\delta_{N2}^I$ and
			(b) $\delta_{ML}$, $\delta_{ML}^I.$}
		\label{fig:1}
	\end{figure}
	\newpage
	\small{
		\hspace{1.3cm}{\bf Table 5.2: $n=8~~~~~~~~~~~~~~~~~~~~~~~~~~~~~$}\vspace{0.05cm}
		\begin{center}
			\begin{tabular}{|c|c||c|c|c||c|c||}
				\hline
				\emph{$\sigma_1^{-1}
					\downarrow$}&\emph{$\sigma_2^{-1}
					\downarrow$}&{\emph{$R(\delta_{N1})$}}&{\emph{$R(\delta_{N2})$}}&{\emph{$R(\delta_{N2}^I)$}}
				&{\emph{$R(\delta_{ML})$}}&{\emph{$R(\delta_{ML}^I)$}}\\
				\hline
				~&0.2&	0.095823& 0.062259&	0.061105&	0.049314&	0.048990\\
				~&0.4&	0.103624& 0.065513&	0.064204&		0.048021&	0.047518\\
				0.3&	0.6& 0.085579&	0.061412&	0.060521&		0.057865&	0.057707\\
				~&0.8&	0.079798& 0.064017&	0.063368&		0.068854&	0.068318\\
				~&1.0& 0.079772&	0.065796&	0.065113&		0.072439&	0.069011\\
				\hline
				~&0.2&	0.081335& 0.06198&	0.061195&		0.063244&	0.063226\\
				~&0.4&	0.102507& 0.063808&	0.062506&		0.045729&	0.045219\\
				0.5&0.6& 0.103574&	0.064113&	0.062787&		0.045271&	0.044733\\
				~&0.8&	0.092039& 0.061425&	0.060376&			0.05143&	0.051227\\
				~&1.0&	0.085209& 0.061156&	0.060274&		0.057722&	0.055718\\
				\hline
				\hline
				~&0.2&		0.08142& 0.067462&	0.06675&	0.074124&	0.074095\\
				~&0.4&	0.087605& 0.059491&	0.058506&		0.051996&	0.051889\\
				0.7&	0.6& 0.102294& 	0.06268&	0.061359&			0.043685&	0.043147\\
				~&0.8& 0.104669& 	0.063982&	0.062600&		0.043914&	0.043321\\
				~&1.0&		0.09512& 0.060274&	0.059096&	0.046047&	0.044638\\
				\hline
				~&0.2&	0.07928& 0.067117&	0.066445&		0.075572&	0.07537\\
				~&0.4&	0.082021& 0.061207&	0.0604&		0.061012&	0.061003\\
				0.9&	0.6& 0.09588&	0.062198&	0.061026&		0.049135&	0.048798\\
				~&0.8& 	0.104904& 	0.063946&	0.062561&	0.043608&	0.043008\\
				~&1.0&	0.106844& 0.065979&	0.064604&		0.045734&	0.043058\\
				\hline
				~&0.2&	0.080529& 0.067564&	0.066761&		0.075219&	0.075167\\
				~&0.4&	0.081238& 0.062603&	0.061841&		0.064586&	0.064401\\
				1.0&	0.6& 0.090518&	0.060502&	0.059452&		0.051105&	0.050918\\
				~&0.8&	0.103641& 0.064743&	0.063423&		0.046463&	0.045936\\
				~&1.0& 0.106738& 	0.065644&	0.064263&		0.045170&	0.044482\\
				\hline
				\hline
			\end{tabular}
		\end{center}
	}
	\begin{figure}[h!]
		\begin{subfigure}{0.6\textwidth}
			\centering
			\includegraphics[width=0.7\linewidth]{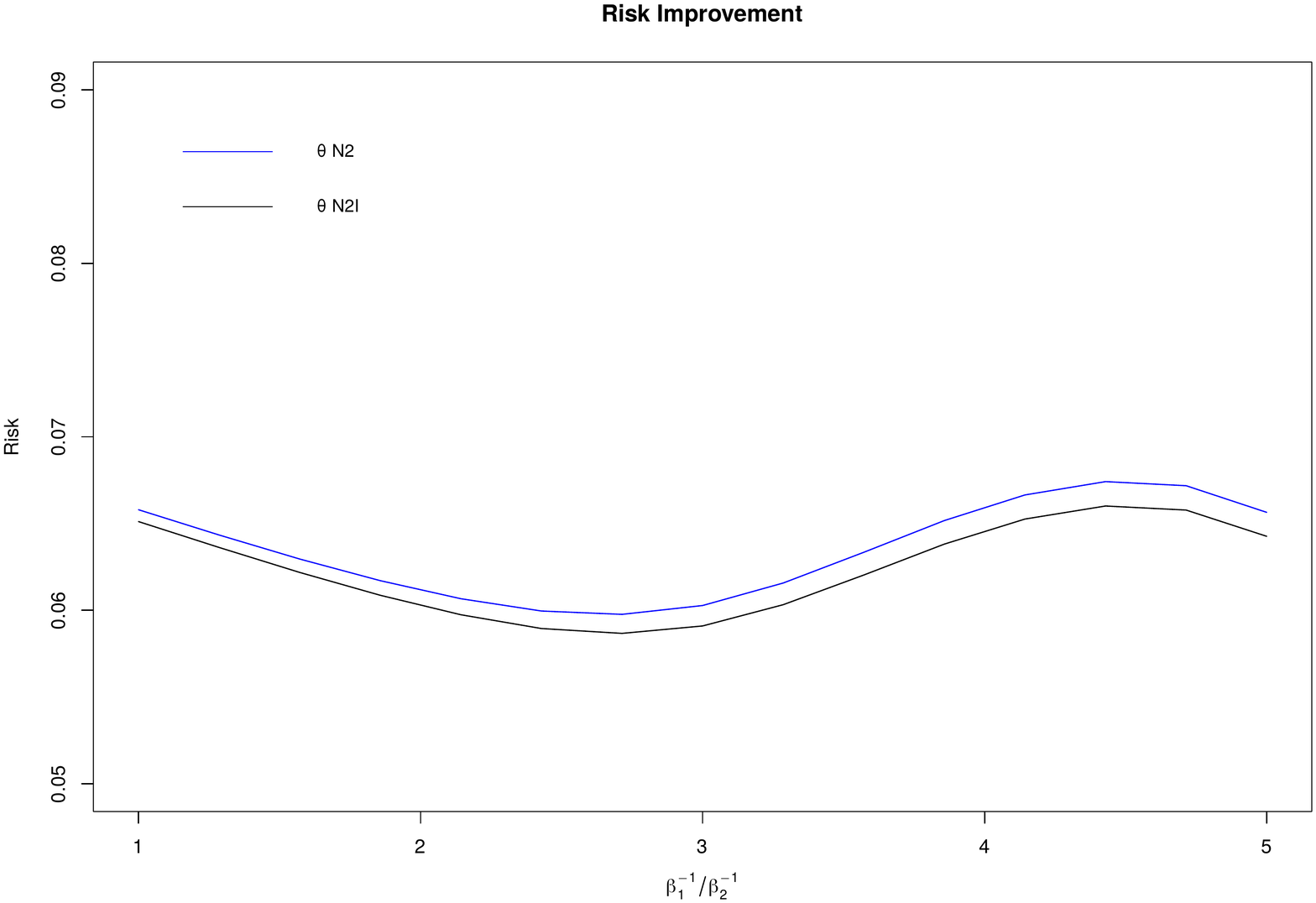}
			\caption{}
			\label{fig:test1}
		\end{subfigure}%
		\begin{subfigure}{0.6\textwidth}
			\includegraphics[width=0.7\linewidth]{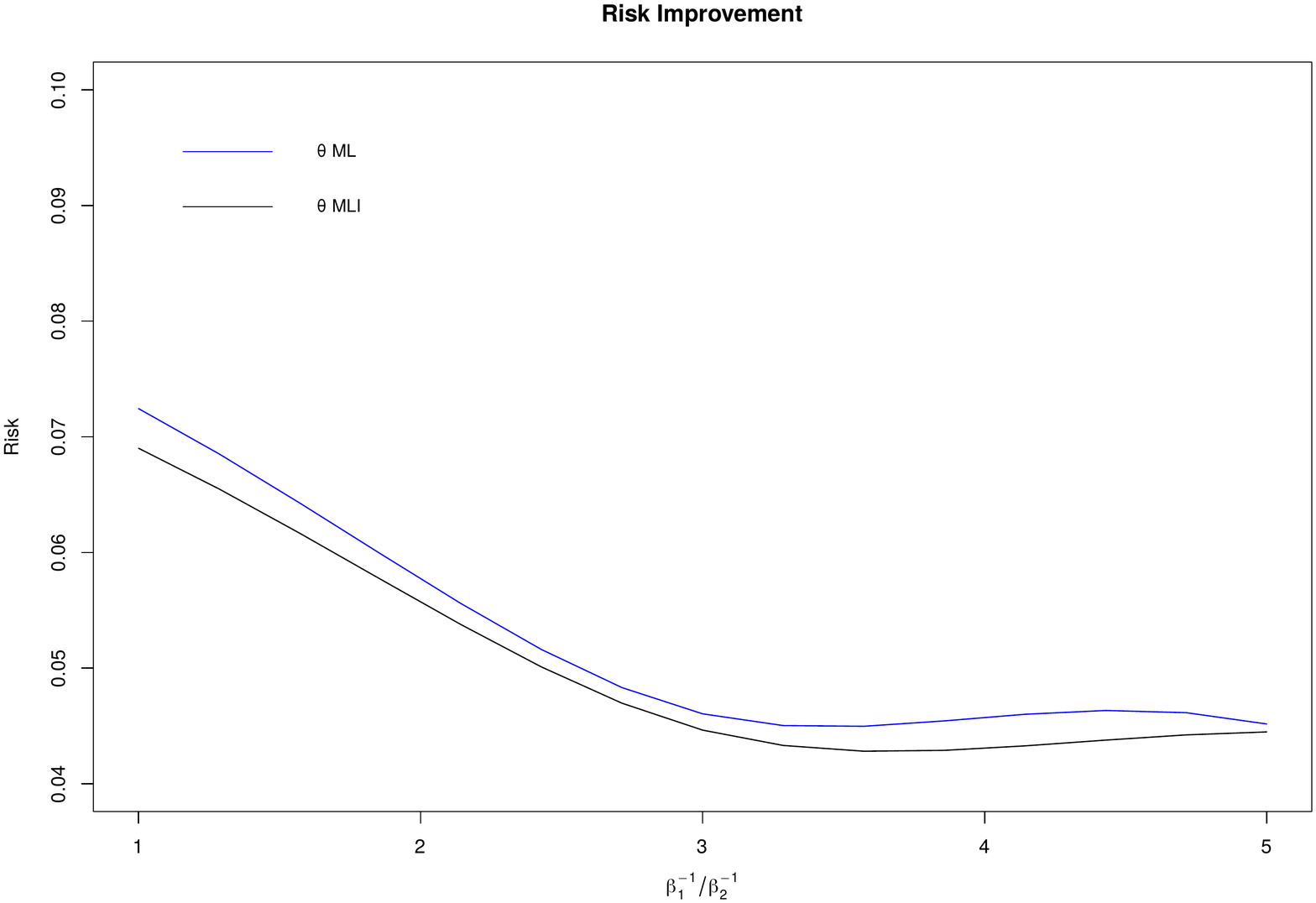}
			\caption{}
			\label{fig:test2}
		\end{subfigure}
		\caption{For $n=8,$ the risk plots of the estimators 
		(a)$\delta_{N2}$, $\delta_{N2}^I$ and
		(b) $\delta_{ML}$, $\delta_{ML}^I.$}
		\label{fig:2}
	\end{figure}
	\section{Conclusion}\label{ch5sec5.6}
	In this article, we have dealt with estimation of hazard rate from a selected exponential population under the entropy loss function. We have obtained the analogue of the best scale equivariant estimator which is proved to be minimax. The inadmissibility of the maximum likelihood estimator and the best scale equivariant estimator has been analyzed. We found an estimator which dominates the minimax estimator $\delta_{N2}$ in terms of risk. It is shown that the estimator $\delta_{N2}^I$ is also minimax. A comparative study for the proposed estimators in terms of risk values has been carried out by numerical simulation. Graphical comparisons of these estimators have been incorporated for better understanding and visualization.
	
	\section*{Acknowledgements}
	The authors are grateful to the editor, associate editor and the anonymous referees for their useful comments and suggestions which led to improvements in the paper.

	\section*{Conflict Of Interest Statement}
	The authors declare that they do not have any conflict of interest.


\begin{thebibliography}{}
		\bibitem{B1} Arshad, M. and Misra, N. (2017): On Estimating the scale parameter
		of the selected uniform population under the entropy loss function. \textit{Brazilian Journal of Probability and Statistics}, {\textbf{31(2)}}, 303-319.
		
		\bibitem{B2} Arshad, M., Misra, N. and Vellaisamy, P. (2014): Estimation after selection
		from Gamma populations with unequal known shape parameters. \textit{Journal of Statistical Theory and Practice}, \textbf{9(2)}, 395-418.
		
		\bibitem{B3} Bahadur, R.R. and Goodman, A.L. (1952): Impartial decision rules and sufficient statistics. \textit{Annals of Mathematical Statistics},
		\textbf{23}, 553-562.
		
		
		\bibitem{B4}Berger, J.O. (1980): Improving on inadmissible estimators in continuous exponential families with applications to simultaneous estimation of gamma scale parameters. \textit{Ann. Statist.}, \textbf{8}, 545-571.
		
		
		\bibitem{B5}Brewster, J.F. and Zidek, J. V. (1974): Improving on equivariant
		estimators. \textit{Ann. Statist.}, \textbf{2}, 21-38.
		
		
		\bibitem{B6} Eaton, M.L. (1967): Some optimum properties of ranking
		procedures. \textit{Ann. Math. Statist.}, \textbf{38}, 124-137.
		
		\bibitem{B7} Jha, B.K., Mahapatra, A.K. and Kayal, S. (2020): Estimation of Hazard Rate of a Selected Exponential Population. \textit{Journal of Statistical Theory and Practice}, \textbf{14(47)}, 1-27.
		
		
		
		\bibitem{B8} Kumar, S., Mahapatra, A.K. and Vellaisamy, P. (2009): Reliability estimation
		of the selected exponential populations. \textit{Statistics and Probability Letters}, \textbf{79}, 1372-1377.
		
		\bibitem{B9} Lehmann, E.L. (1966): On a theorem of Bahadur and Goodman.
		\textit{Ann. Math. Statist.}, \textbf{37}, 1-6.
		
		
		
		\bibitem{B10} Mahapatra, A.K., Kumar, S. and Vellaisamy, P. (2012): Simultaneous estimation of hazard rates of several exponential populations. \textit{Statistica Neerlandica} \textbf{66(2)}, 121–132.
		
		\bibitem{B11} Misra, N., Van der Meulen, E. C. and Vanden Braden, K. (2006a): On estimating the scale parameter of the selected gamma population under the scale invariant squared error loss function. \textit{Journal of Computational and Applied Mathematics}, \textbf{186}, 268-282.
		
		\bibitem{B12} Misra, N., Van der Meulen, E. C. and Vanden Braden, K. (2006b): On some inadmissibility results for the scale
		parameters of selected gamma populations. \textit{Journal of Statistical Planning and
			Inference}, \textbf{136}, 2340 – 2351.
		
		\bibitem{B13} Nematollahi, N. and Motamed-Shariati, F. (2009): Estimation of the Scale Parameter of the Selected
		Gamma Population Under the entropy loss function. \textit{Communications in Statistics - Theory and Methods}, \textbf{38(2)}, 208-221.
		
		\bibitem{B14} Parsian, A. and Nematollahi, N. (1996):Estimation of scale parameter under entropy
		loss function, \textit{Journal of Statistical Planning and
			Inference}, \textbf{52}, 77-91.
		
		
		
		\bibitem{B15} Sackrowitz, H.B. and Samuel-Cahn, E. (1986): Evaluating the chosen
		population: A Bayes and Minimax approach. \textit{Adaptive Statistical
			Procedures and Related topics, IMS Lecture Notes-Monograph Series}, \textbf{8}, 386-399.
		
		\bibitem{B16} Sharma, D. (1977): Estimation of the reciprocal of the scale parameter in a shifted exponential
		distribution. \textit{Sankhya-A}, \textbf{39(2)}, 203-205. 
		
		
		
		
		\bibitem{B17} Vellaisamy, P. (1992): Inadmissibility results for the selected
		scale parameters. \textit{Ann. Statist.}, \textbf{20}, 2183-2191.
		
		
		\bibitem{B18} Vellaisamy, P. (1996): A note on the estimation of the selected
		scale parameters. \textit{J. Statist. Plann. Inference}, \textbf{55}, 39-46.
		
		
		
	\end{thebibliography}
\end{document}